\documentclass[graybox]{svmult}

\usepackage{mathptmx}       
\usepackage{helvet}         
\usepackage{courier}        
\usepackage{type1cm}        

\usepackage{makeidx}         
\usepackage{graphicx}        
\usepackage{multicol}        
\usepackage[bottom]{footmisc}
\usepackage{amssymb}
\usepackage{amsfonts}
\usepackage{amsmath}
\usepackage{textcomp}
\usepackage{xfrac}
\usepackage{subcaption}

\newcommand{\rP}{\mathrm{P}} 
\newcommand{\rE}{\mathrm{E}} 

\makeindex             

\begin{document}

\title{A multitype Galton-Watson model for rejuvenating  cells
}

\author{Serik Sagitov, Lotta Eriksson and Marija Cvijovic}
\institute{Serik Sagitov \at Department of Mathematical Sciences, Chalmers University of Technology and University of Gothenburg SE-412 96 Göteborg
Sweden, \email{serik@chalmers.se}
\and Lotta Eriksson \at Department of Mathematical Sciences, Chalmers University of Technology and University of Gothenburg SE-412 96 Göteborg
Sweden, \email{lottaer@chalmers.se}
\and Marija Cvijovic \at Department of Mathematical Sciences, Chalmers University of Technology and University of Gothenburg SE-412 96 Göteborg
Sweden, \email{marija.cvijovic@chalmers.se}
}
%
%
\maketitle

\abstract*{We employ the framework of multitype Galton-Watson processes to model a population of dividing cells. The cellular type is represented by its biological age, defined as the count of harmful proteins hosted by the cell. The stochastic evolution of the biological age of a cell is modeled as a discrete Markov chain with a finite state space $\{0,1,\ldots,n\}$, where $n$ signifies the absorbing state corresponding to the senescent state of the cell. 
Consequently, the set of individual types in the multitype Galton-Watson process becomes $\{0,\ldots,n-1\}$.
In our setting, a dividing cell may undergo rejuvenation meaning that its biological ages reduces due to transfer of harmful proteins to the daughter cell. For the proposed model, we define and study several biologically meaningful features, such as rejuvenation states, expected replicative lifespan, population growth rate, stable biological age distribution, and the population average of the biological age.}

\abstract{We employ the framework of multitype Galton-Watson processes to model a population of dividing cells. The cellular type is represented by its biological age, defined as the count of harmful proteins hosted by the cell. The stochastic evolution of the biological age of a cell is modeled as a discrete Markov chain with a finite state space $\{0,1,\ldots,n\}$, where $n$ signifies the absorbing state corresponding to the senescent state of the cell. 
Consequently, the set of individual types in the multitype Galton-Watson process becomes $\{0,\ldots,n-1\}$. \\
In our setting, a dividing cell may undergo rejuvenation, meaning that its biological age reduces due to the transfer of harmful proteins to the daughter cell. For the proposed model, we define and study several biologically meaningful features, such as rejuvenation states, expected replicative lifespan, population growth rate, stable biological age distribution, and the population average of the biological age.}

\section{Introduction}

Aging is a complex biological phenomenon driven by the progressive accumulation of cellular damage, which ultimately manifests as declining physiological functions and increased susceptibility to age-associated diseases. The gradual accumulation of damage disrupts essential biomolecules, thereby compromising the functionality of both individual cells and the entire organism. It is frequently seen as both the core and driving force of aging.
Accumulated damage in the parent cells contributes to initial damage at birth, burdening individuals with reduced lifespan and health span. In unicellular systems, such as budding yeast \textit{Saccharomyces cerevisiae}, damaged proteins exhibit asymmetric inheritance, primarily retained within the mother cell through active transport mechanisms such as retention \cite{Hugo}. During cell division, this process selectively segregates damaged cellular components within the mother cell compartment, which is
crucial for maintaining cell population viability and generating rejuvenated daughter cells \cite{Kennedy, Lippuner}. Rejuvenation, thus, can be seen as a reduction in the amount of inherited damage in successive generations, practically making them more youthful. This, in a broader sense, can imply that the measurable amount of accumulated damage is a proxy for biological age. Previously,  aging phenomena have been studied via dynamic and multi-scale modeling, spanning single-cell and population levels, to elucidate the accumulation of damaged proteins and its impact on fitness, rejuvenation, health span, and aging within asymmetrically \cite{Qasim, JohannesNiek, SBC, BarbaraLinnea} and symmetrically \cite{Ari, Ack1, Ack2} dividing yeast cells.

In this work, we propose an individual-based stochastic population model for cells whose probabilities of division depend on their biological ages. The biological age of a cell is assumed to evolve according to a time-homogeneous discrete Markov chain $\xi=(\xi_0,\xi_1,\ldots)$ with finitely many states $\{0,1,\ldots,n\}$, introduced in Section \ref{rw}. The state  $\xi_t=i$  represents the current biological age of the cell, assumed as the number of harmful proteins within the cell at time $t$. The state $n$ is an absorption state corresponding to the critical number of harmful proteins at which the cell is unable to divide. At any other state $i\in\{0,\ldots,n-1\}$, the cell receives a random number $\tau$ of additional harmful proteins and then divides with a certain probability $b_{i+\tau}$.
During the birth event, the mother cell releases a fraction $1-p$ of harmful proteins to its daughter cell, potentially leading to rejuvenation of the mother cell $\xi_{t+1}<\xi_t$.  Accordingly, a Markov chain state $i$ will be called a \textit{rejuvenation state} if 
$$\rE(\xi_1|\xi_0=i)< i,$$
meaning that for the cell of biological age $i$ at some time $t$, the \textit{expected} biological age at the time $t+1$ is less then $i$. Proposition \ref{MC} of Section \ref{rw} identifies the set of such rejuvenation states in terms of 
\begin{itemize}
    \item the senescence age, denoted as $n$,
    \item the inflow per cell of harmful proteins, denoted as $\tau$,
    \item the retention proportion, denoted as $p$,
    \item the sequence of division probabilities, denoted as $(b_0,b_1,\ldots)$,
\end{itemize}
where 
$n$ is a positive integer number, $\tau$ is a non-negative integer valued random variable,  and $p\in[0,1]$.

In Section \ref{mal}, we present the key special case, which we call the $(m,\alpha)$-case, that allows for more detailed analysis and simulation studies. It is described by four parameters $(n,m,p,\alpha)$, where $m$ is a constant inflow and the parameter $\alpha\ge1$ fully determines the sequence of division probabilities $(b_i)$.

In Section \ref{rL}, we specify the replicative lifespan of a cell as a pertinent functional associated with the Markov chain $\xi$. The replicative lifespan of a cell is defined as the total number of divisions the cell undergoes before reaching the state of senescence \cite{Steinkraus}. For a fixed $n$ and the current value of the biological age $i$, we study the expected replicative lifespan $\lambda_i$ as a function of the model parameters $(m,p,\alpha)$. 

\begin{figure}[h]
\centering
 \includegraphics[width=\textwidth]{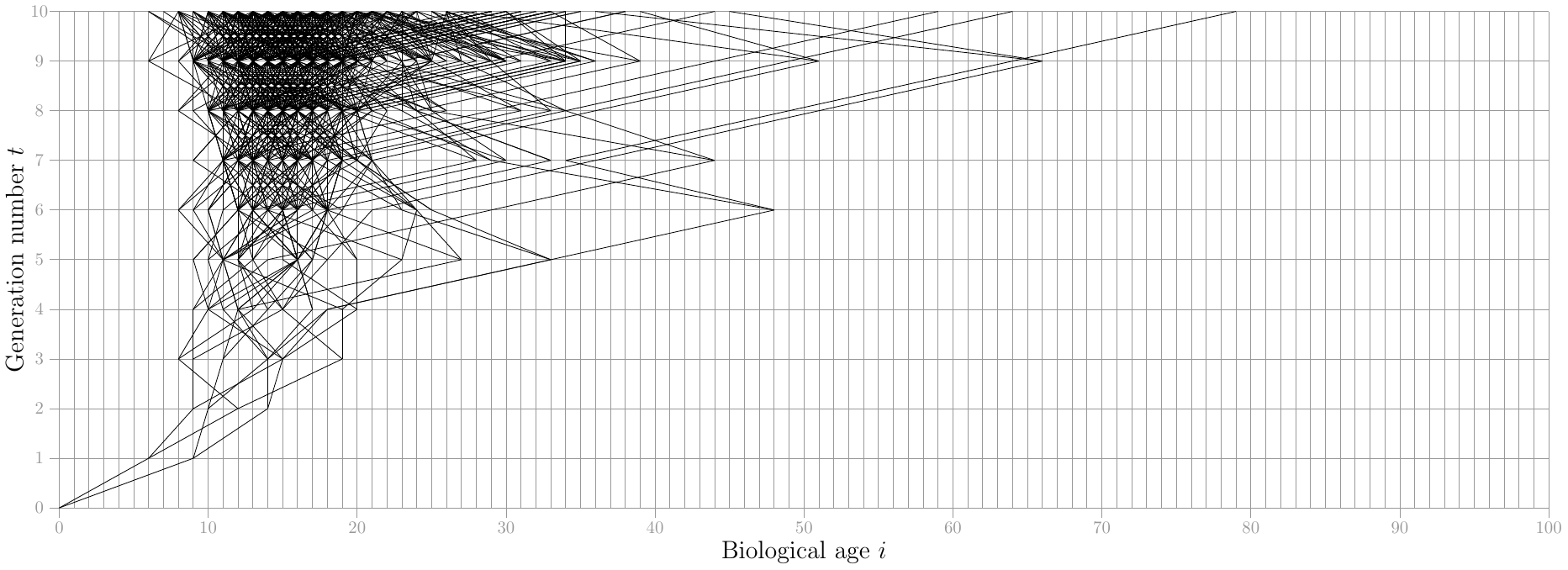}
 \caption{A tree realization for the $(n,m,p,\alpha)$-model  with parameters $n=100$, $m=15$, $p=0.5$, and $\alpha=2$. 
  }
 \label{ill}
\end{figure}

The stochastic population model is presented in Section \ref{zz} as a multitype Galton-Watson process \cite{Hac, J, Kim}, with $n$ types of individuals. Here, a Galton-Watson individual of type $i$ represents a non-senescent cell hosting  $i\in \{0,\ldots,n-1\}$ harmful proteins. It is crucial to emphasize that our approach relies on the labelling of individual types with biological age, deviating from the utilization of either chronological age, as in \cite{JS}, or the number of telomeres, as in \cite{OK}. The multitype Galton-Watson process is a vector-valued Markov chain
\begin{equation}\label{Zt}
\boldsymbol{Z}_t=(Z_t(0),\ldots,Z_t(n-1)), \quad t=0,1,\ldots, 
\end{equation}
where $Z_t(i)$ stands for the number of cells of biological age $i$ in the population at time $t$.  A tree representation of a trajectory of such a stochastic reproduction process showing an outcome for ten generations simulated using the parameter values $(n,m,p,\alpha)=(100, 15,0.5,2)$ is depicted in Figure \ref{ill}. 

In Section \ref{PAA}, we discuss the average biological age across a population exhibiting a stable distribution of the cell types. We reveal different patterns of the relationship between the average biological age and the growth rate of a cell population depending on the parameter set $(n,m,p,\alpha)$. Section \ref{proof} includes various proofs.


\section{Rejuvenation states of the Markov chain $\xi$}\label{rw}
Suppose we are given the values of the senescent state $n$, the discrete distribution
\[\rP(\tau=i)=q_i, \quad q_0+q_1+\ldots=1,\] 
of the random inflow $\tau\ge0$, the probability of retention $0\le p\le1$, and a sequence of division probabilities $(b_i)$. In this section, we define a Markov chain $\xi$, describing the evolution of the biological age of a single mother cell, together with an accompanying sequence 
\[\eta=(\eta_1,\eta_2,\ldots),\]
describing the types of eventual daughter cells. We will write $i\wedge n$ instead of $\min(i,n)$ and $B_k$ will mean a realisation of a binomial random variable Bin$(k,p)$.

The transition rule of the Markov chain $\xi$ having the state space $\{0,\ldots,n\}$, represents the following chain of events during the time interval $[t,t+1]$:
\begin{itemize}
    \item given $\xi_t=i$ with $0\le i<n$, after the mother cell receives $\tau=k$ additional harmful proteins,  it divides with  probability  $b_{i+k}$,
    \item  in the case of no division, tagged by the relation $\eta_{t+1}=n+1$, the mother cell becomes a host of $\xi_{t+1}=(i+k)\wedge n$ harmful proteins,
    \item  in the case of division, the mother cell becomes a host of $\xi_{t+1}=B_{i+k}\wedge n$ harmful proteins, and the daughter cell is born with $\eta_{t+1}=(i+k-B_{i+k})\wedge n$ harmful proteins,
    \item given $\xi_t=n$, the cell  stops evolving resulting in $\xi_{t+1}=n$ and $\eta_{t+1}=n+1$.
\end{itemize}
Here, the binomial distribution appears due to the assumption that in the case of cell division, each of the $i+k$ harmful proteins is either retained by the mother cell with probability $p$ or transferred to the daughter cell with probability $1-p$.
In summary, the state $n$ acts as an absorption state of the Markov chain $\xi$, and if $\xi_t=i$ with $0\le i<n$, and $\tau=k$, then
\begin{equation}\label{trans}
\{\xi_{t+1};\eta_{t+1}\}=
\left\{
\begin{array}{lll}
 \{(i+k)\wedge n;n+1\} &   \text{with probability }& 1-b_{i+k},   \\
 \{B_{i+k}\wedge n;(i+k-B_{i+k})\wedge n\} &   \text{with probability }& b_{i+k},\end{array}
\right.
\end{equation} 
where 
in the case 
$(i+\tau)\wedge n=(B_{i+\tau})\wedge n=j$,
the rule \eqref{trans} implies $\rP(\xi_{t+1}=j)=1$. 

Denote by
\[\chi_{t}=1_{\{\eta_t\le n\}},\]
the indicator of the cell division outcome and turn to 
the transition probabilities  of the Markov chain $\xi$
 \[p_{ij}=\rP(\xi_{t+1}=j|\xi_t=i),\quad 0\le i,j\le n.\]
Putting
\begin{align*}
    Q_{ij}&=\rP(\xi_{t+1}=j,\chi_{t+1}=0|\xi_t=i),\\
    R_{ij}(p)&=\rP(\xi_{t+1}=j,\chi_{t+1}=1|\xi_t=i),
\end{align*}
we obtain
\begin{equation}\label{pij}
p_{ij}=Q_{ij}+R_{ij}(p).
\end{equation}
Let
\[g_{kj}(p)={k\choose j}p^j(1-p)^{k-j}1_{\{0\le j\le k\}}.\]

\begin{proposition}\label{MC}
 Relation \eqref{pij} holds for $i\le j<n$, with
\[Q_{ij}=q_{j-i}(1-b_j),\quad R_{ij}(p)=\sum_{k\ge j}q_{k-i}b_kg_{kj}(p),\]
for $0\le j<i<n$, with
\[Q_{ij}=0,\quad R_{ij}(p)=\sum_{k\ge i}q_{k-i}b_{k}g_{kj}(p),\]
and for $0\le i<j=n$, with 
\[Q_{in}=\sum_{k\ge n}q_{k-i}(1-b_k),\quad  R_{in}(p)=\sum_{k\ge n}q_{k-i}b_k\sum_{j=n}^{k}g_{kj}(p).\]
\end{proposition}

Set
\[h_i=\sum_{j=1}^{n}jp_{ij}.\]
The probability of a rejuvenation event  $\xi_{t+1}<\xi_t$ depends, of course, on the current state of Markov chain $\xi$. Therefore, it is relevant to introduce the concept of a rejuvenation state as follows: the state $i$ will be called a rejuvenation state of the Markov chain $\xi$ if $d_i=h_i-i$, the expected jump size from the state $i$ is negative.
Clearly, $d_0>0$ and $d_n=0$, implying that the set of the rejuvenation states belongs to $\{1,\ldots,n-1\}$.  

\begin{proposition}\label{coro}
Consider the Markov chain $\xi$ assuming  $b_i=0$ for $i\ge n$. Then 
\begin{equation}\label{dav}
h_i=\rE((i+\tau)\wedge n -(1-p)(i+\tau)b_{i+\tau}),
\end{equation}
and the set of rejuvenation states is defined by the relations
\begin{equation}\label{day}
\rE(\tau\wedge (n-i))< (1-p)\rE((i+\tau)b_{i+\tau}),\quad 1\le i\le n-1.
\end{equation}
\end{proposition}

\section{$\boldsymbol{(m,\alpha)}$-case}\label{mal}

\begin{figure}[h!]
\centering
 \includegraphics[width=8cm]{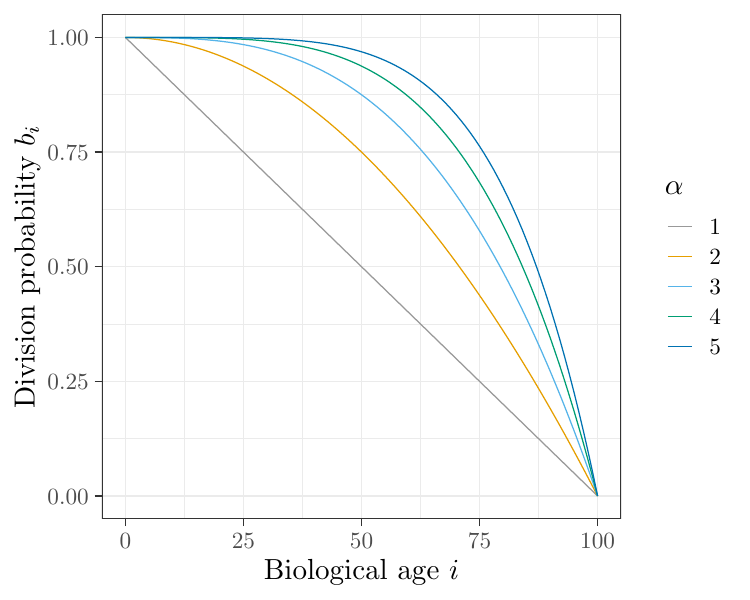}
 \caption{An illustration of  formula \eqref{spe} with $n=100$. The five curves, arranged in ascending order, depict five distinct division rates $\alpha= 1, 2, 3, 4,5$. }
\label{bill}
\end{figure}

In this paper, special attention is paid to the case when the division probabilities satisfy 
\begin{equation}\label{spe}
    b_{i}=(1-(i/n)^\alpha)1_{\{0\le i\le n\}},\quad \alpha\ge1.
\end{equation} 
The defining parameter $\alpha$ will be called the division rate, justified by the fact that for larger values of $\alpha$, the probabilities of cell division become uniformly larger, as illustrated by Figure \ref{bill}. The case, when besides \eqref{spe}, it is also assumed that the inflow $\tau$ is a constant $m$,
\begin{equation}\label{spm}
    \rP(\tau=m)=1,\quad m\in\{1,\ldots,n-1\},
\end{equation} 
will be called the $(m,\alpha)$-case. This case lends itself well to studying both analytically and through simulations.

Observe that in the $(m,\alpha)$-case, for $i\ge n-m$, we have $b_{i+m}=0$ and
\[p_{in}=1,\quad  n-m\le i\le n.\]
On the other hand, for $0\le i<n-m$, relation \eqref{pij} holds with
\begin{equation}\label{nel}
Q_{ij}=(\tfrac{i+m}{n})^\alpha1_{\{j=i+m\}},\quad R_{ij}(p)=(1-(\tfrac{i+m}{n})^\alpha)g_{i+m,j}(p).
\end{equation} 


\begin{proposition}\label{prog} Consider the Markov chain $\xi$ assuming \eqref{spe} and \eqref{spm}. Then the expected jump size is  $d_i=n-i$ for $i=n-m,\ldots,n$, and
\[d_i=m-(1-p)(i+m)(1-(\tfrac{i+m}{n})^\alpha),\quad 0\le i<n-m.\]
In this $(m,\alpha)$-case, the set of rejuvenation states is either empty or forms an interval $i_1\le i\le i_2$. In the latter case, the end of the intervals are such that $0< i_1\le i_2<n-m$ and 
\[i_1=\lceil ny_1\rceil-m,\quad i_2=\lfloor ny_2\rfloor-m,\]
where $y_1< y_2$ are two distinct roots of the equation
\[y-y^{1+\alpha}= \tfrac{m}{(1-p)n},\quad 0<y<1.\]
\end{proposition}

 \subsubsection*{Example} 
 To illustrate Proposition \ref{prog}, we turn to the set of parameters used for generating Figure \ref{ill}: 
 $$n=100, \quad m=15, \quad p=0.5,\quad \alpha=2.$$
In this particular case, after solving the equation
\[y-y^3= 0.3\quad 0<y<1,\]
we identify the set of rejuvenation states as $19\le i\le 63$.

\begin{proposition}\label{cop}
Consider $i_1=i_1(p,m,\alpha)$ and $i_2=i_2(p,m,\alpha)$ from Proposition \ref{prog} as functions of three variables for the Markov chain $\xi$ assuming \eqref{spe} and \eqref{spm}. Then the following monotonicity properties hold
\begin{align*}
 i_1(p,m,\alpha)&<  i_1(p+\epsilon,m,\alpha)< i_2(p+\epsilon,m,\alpha)<  i_2(p,m,\alpha),\\
 i_1(p,m,\alpha)&< i_1(p,m+1,\alpha)< i_2(p,m+1,\alpha)< i_2(p,m,\alpha),\\
 i_1(p,m,\alpha+\epsilon)&< i_1(p,m,\alpha)< i_2(p,m,\alpha)<  i_2(p,m,\alpha+\epsilon),
\end{align*}
where $\epsilon>0$ and the corresponding values of $i_1, i_2$ are well defined.
\end{proposition}

\begin{figure}[h!]
    \centering
    \includegraphics[width=0.7\textwidth]{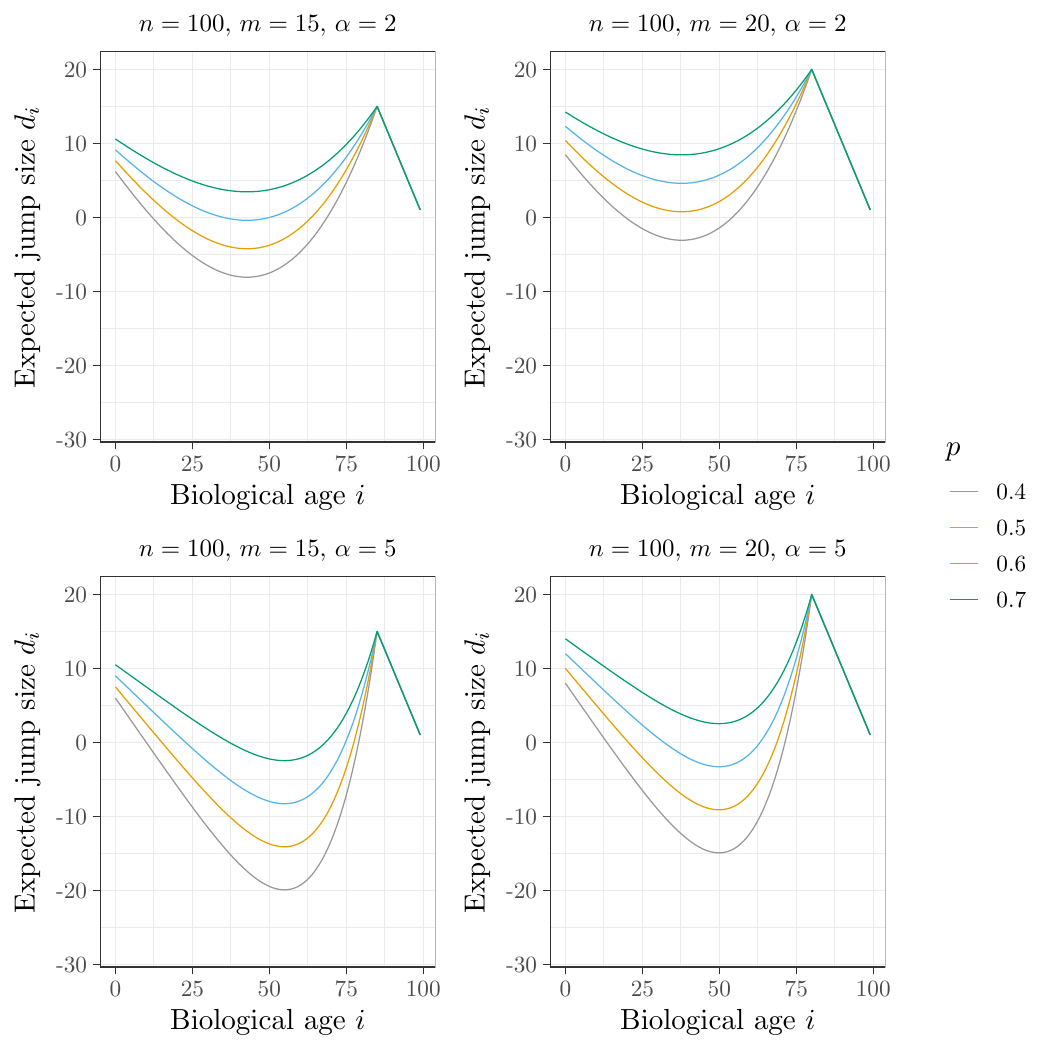}
    \caption{The expected jump size $d_i$ is depicted as a function of the biological age $i$ assuming $n=100$, $m=15$ or $m=20$, $\alpha=2$ or $\alpha=5$. The four curves in each of the panels, arranged in ascending order, represent four values of the retention proportion $p=0.4, 0.5, 0.6, 0.7$.}
    \label{fdi}
\end{figure}

Figure \ref{fdi} illustrates Proposition \ref{cop}. For instance, in the top-left panel, an increase in $p$ elevates the curve $d_i$, resulting in the contraction of the interval associated with the inequality $d_i<0$. A similar trend is observed when comparing the panels on the left and right, as the inflow parameter shifts from $m=15$ to $m=20$.

\section{Expected replicative  lifespan}\label{rL}
The replicative lifespan $L$, interpreted as the total number of cell divisions until its senescence, can be defined  as the following functional of the Markov chain $\xi$:
\[L=\sum_{t\ge1}\chi_t.\]
By construction, the random variable $L$ has a phase-type distribution \cite[Ch. 3]{As}.
Turning to the expected replicative lifespan of a mother cell of age $i$ 
$$\lambda_i=\rE(L|\xi_0=i),$$ 
observe first that  $\lambda_n=0$, since $n$ is the absorption state of the chain $\xi$. The Markov property of the chain $\xi$ entails the recursion
\begin{equation}
 \label{rig}
\lambda_i=\beta_i+\sum_{j=0}^{n-1} p_{ij}\lambda_j,\quad i=0,\ldots,n-1,
\end{equation}
where $\beta_i=\rE(b_{i+\tau})$.
Letting
\[\boldsymbol{\lambda}=(\lambda_0,\ldots,\lambda_{n-1}),\quad \boldsymbol{\beta}=(\beta_0,\ldots,\beta_{n-1}),\quad \mathbb P_0=(p_{ij})_{i,j=0}^{n-1},\]
we can rewrite recursion \eqref{rig} in the matrix form
$\boldsymbol{\lambda}^\intercal=\boldsymbol{\beta}^\intercal+\mathbb P_0 \, \boldsymbol{\lambda}^\intercal$
yielding
\begin{equation}
 \label{rg}
\boldsymbol{\lambda}^\intercal=(\mathbb I-\mathbb P_0)^{-1}\boldsymbol{\beta}^\intercal.
\end{equation}

\subsection*{$\boldsymbol{(m,\alpha)}$-case}
Given \eqref{spe} and \eqref{spm}, we get $\lambda_{n-m}=\ldots=\lambda_{n-1}=0$ and
\[\lambda_i=b_{i+m}+\sum_{j=0}^{n_1} p_{ij}\lambda_j,\quad i=0,\ldots,n_1,\]
where $n_1=n-m-1$. This leads to a truncated version of \eqref{rg},
\begin{equation}
 \label{rog}
\boldsymbol{\hat \lambda}^\intercal=(\mathbb I-\mathbb {\hat P}_0)^{-1} \boldsymbol{\hat\beta}^\intercal,
\end{equation}
with
\[\boldsymbol{\hat \lambda}=(\lambda_0,\ldots,\lambda_{n_1}),\quad \boldsymbol{\hat \beta}=(b_m,\ldots,b_{n-1}),\quad \mathbb 
{\hat P}_0=(p_{ij})_{i,j=0}^{n_1},\]
where $p_{ij}$ are given by \eqref{pij} with \eqref{nel}.

\begin{figure}[h!]
    \centering
    \includegraphics[width=0.75\textwidth]{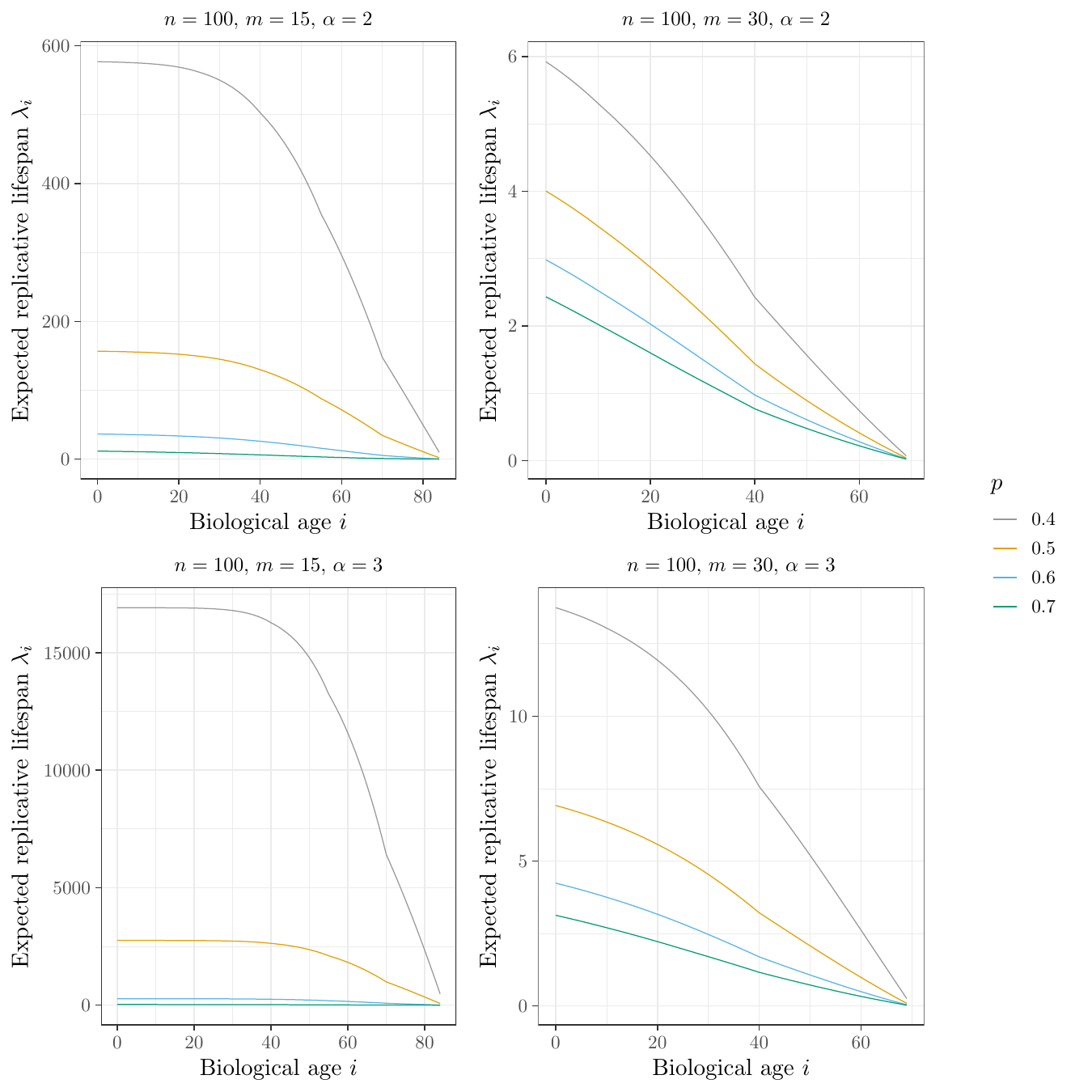}
    \caption{The expected replicative lifespan $\lambda_i$  assuming $n=100$, $m=15,30$, $\alpha=2,3$. The four curves in each of the panels, arranged in descending order, represent $p=0.4, 0.5, 0.6, 0.7$.}
    \label{flam}
\end{figure}

\begin{proposition}\label{coup}
For a fixed $n$, consider the expected replicative lifespan $\lambda_i=\lambda_i(p,m,\alpha)$ as a function of four variables assuming \eqref{spe} and \eqref{spm}. For each of the variables this function is monotone in that 
\begin{align*}
     \lambda_i(p,m,\alpha)&> \lambda_{i+1}(p,m,\alpha),\quad \lambda_i(p,m,\alpha)> \lambda_i(p+\epsilon,m,\alpha),\\ 
     \lambda_i(p,m,\alpha)&> \lambda_i(p,m+1,\alpha),\quad
     \lambda_i(p,m,\alpha)< \lambda_i(p,m,\alpha+\epsilon),
\end{align*}
whenever $0\le i<n-m$, $0\le p<p+\epsilon\le1$, and $1\le \alpha<\alpha+\epsilon$.
\end{proposition}

The assertions of Proposition \ref{coup}, illustrated by Figure \ref{flam},  are quite intuitive: the expected replicative lifespan increases if we either decrease the current biological age $i$, the retention proportion $p$, and the inflow $m$ or if we increase the division rate $\alpha$. 

Examining the top left panel of Figure \ref{flam}, we observe that the expected replicative lifespan for parameters $n=100$, $m=15$, $p=0.5$, and $\alpha=2$ typically falls within the range of 100 to 150. Returning to Figure \ref{ill}, we can conclude that observing the population for just 10 time units is insufficient to witness cells entering the senescent state.

\section{Multitype Galton-Watson process}\label{zz}

Taking the transition rule outlined in Section \ref{rw} for the Markov chain $\xi$ and the associated sequence $\eta$ as a representation of the evolutionary dynamics for $n$ distinct types of cells, and assuming these cells evolve independently of each other, we obtain the multitype Galton-Watson process \eqref{Zt}. The vector of expected population sizes at time $t$ 
\[\rE(\boldsymbol{Z}_t|\boldsymbol{Z}_0)=\mathbb M^t\boldsymbol{Z}_0^\intercal,\]
is expressed in terms of the matrix
\begin{equation}
 \label{ma}
\mathbb M=(M_{ij})_{i,j=0}^{n-1},
\end{equation}
consisting of  the mean numbers  of offspring of type $j$ produced by a mother of type $i$
\begin{align*}
 M_{ij}&=\rE(1_{\{\xi_1=j\}}+1_{\{\eta_1=j\}}|\xi_0=i)=p_{ij}+\rP(\eta_1=j|\xi_0=i).
\end{align*}
Recalling \eqref{pij}, we obtain
\begin{align*}
 M_{ij}&=Q_{ij}+R_{ij}(p)+R_{ij}(1-p),
\end{align*}
implying that the mean offspring numbers symmetrically depend on the parameter $p$ around the central value of $p=0.5$.

The growth rate of the Galton-Watson process \eqref{Zt} is determined by the Perron-Frobenius eigenvalue $r$ of the mean matrix $\mathbb M$, in that
\begin{equation}
 \label{tuv}
\mathbb M^t\sim r^t \boldsymbol{u}^\intercal \boldsymbol{v},\quad t\to\infty,
\end{equation}
where $\boldsymbol{v}=(v_0,\ldots, v_{n-1})$ and $\boldsymbol{u}=(u_0,\ldots, u_{n-1})$ are left and right  eigenvectors such that
\[\mathbb M \boldsymbol{u}^\intercal=r \boldsymbol{u}^\intercal,\quad \boldsymbol{v}\mathbb M=r \boldsymbol{v},\quad \boldsymbol{v}\boldsymbol{u}^\intercal=1.\]

\subsection*{$\boldsymbol{(m,\alpha)}$-case}\label{erw}

Turning to \eqref{nel} and recalling notation $n_1=n-m-1$, we find that 
\begin{equation}
 \label{muv}
M_{ij}=(\tfrac{i+m}{n})^\alpha1_{\{j=i+m\}}+(1-(\tfrac{i+m}{n})^\alpha)(g_{i+m,j}(p)+g_{i+m,j}(1-p))
\end{equation}
for $0\le i\le n_1$, $0\le j\le i+m$, and $M_{ij}=0$ otherwise. 
\begin{proposition}\label{linn}
    Given \eqref{spe}, \eqref{spm}, and $0<p<1$, consider two matrices
 $$\mathbb{\hat M}=(M_{ij}: 0\le i,j\le n_1),\quad \mathbb{M}^*=(M_{ij}: 0\le i<n-m\le j<n),$$ 
 where $M_{ij}$ are specified by \eqref{muv}.
Denote by $r$ the Perron-Frobenius eigenvalue of the square matrix $\mathbb{\hat M}$. Let  $\boldsymbol{\hat v}=(v_0,\ldots, v_{n_1})$ and $\boldsymbol{\hat u}=(u_0,\ldots, u_{n_1})$ be such that
\[\mathbb {\hat M} \boldsymbol{\hat u}^\intercal=r \boldsymbol{\hat u}^\intercal,\quad \boldsymbol{\hat v}\mathbb {\hat M}=r \boldsymbol{\hat v},\quad \boldsymbol{\hat v}\boldsymbol{\hat u}^\intercal=1.\]    
Then \eqref{tuv} holds with
\[\boldsymbol{u}=(u_0,\ldots, u_{n_1},0,\ldots,0),\quad \boldsymbol{v}=(v_0,\ldots,v_{n-1}),\quad (v_{n-m},\ldots,v_{n-1})=r^{-1}\boldsymbol{\hat v}\mathbb{M^*}.\]
    \end{proposition}

 \subsubsection*{Example} 
Returning to Figure  \ref{ill}, we point out that in the case $(n,m,p,\alpha)=(100,15,0.5,2)$, the corresponding interval of the rejuvenation states $19\le i\le 63$ acts as a potential well within the type space. We see that during the first 10 generations of the reproduction process, no cell was able to escape the potential well. Such an escape may require no-division outcomes in three consecutive time units for the same mother cell.


\begin{figure}[h!]
\centering
 \includegraphics[width=\textwidth]{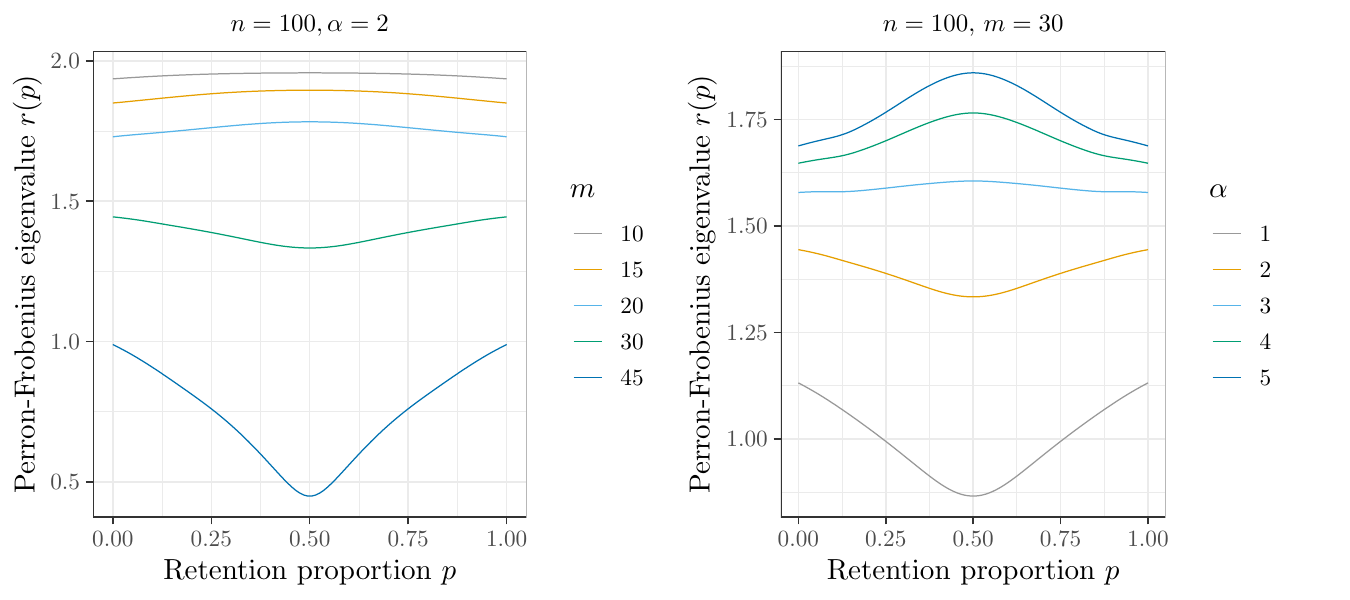}
 \caption{The Perron-Frobenius eigenvalue $r=r(p)$ as the function of the retention proportion $p\in(0,1)$.  Left panel: the five curves are generated for $n=100$, $\alpha=2$, and   $m= 10, 15, 20, 30, 45$ in the descending order of the curves. Right panel: the five curves are generated for $n=100$, $m= 30$, and  $\alpha=1, 2, 3, 4, 5$ in the ascending order of the curves.
  }
 \label{fig3}
\end{figure}

Figure \ref{fig3} depicts the curves $r=r(p)$ for various combinations of the model parameters $(m,p,\alpha)$, illustrating Proposition \ref{linn}. The left panel demonstrates that the growth rate increases if we decrease the inflow $m$ of the harmful proteins. The right panel shows that the growth rate increases with the division rate $\alpha$. As expected, all curves are symmetrical over $p$ around the middle value $p=0.5$. Notice the upper bound $r<2$ evident in both panels. This limitation is easily explained by the model assumption that each Galton-Watson individual produces at most two offspring.

The lowest curve in the left panel of Figure \ref{fig3} corresponding to the parameter set $n=100$, $m=45$, and $\alpha=2$, presents the subcritical case when $r<1$ and the Galton-Watson population goes extinct with probability one. The subcritical regime can be found even in the case  $n=100$, $m=30$, and $\alpha=1$ for the values of $p$ close to 0.5, as can be seen from the lowest curve on the right panel of Figure \ref{fig3}.

Another noteworthy aspect of the curves depicted in Figure \ref{fig3} is that various combinations of the model parameters $(n, m, \alpha)$ can result in the symmetric allocation ($p=0.5$) leading to either the highest or lowest growth rate of the population.

\begin{proposition}\label{cou}
Given \eqref{spe}, \eqref{spm}, for a fixed pair $(n,p)$ with $0<p<1$, consider the Perron-Frobenius growth rate $r=r(m,\alpha)$ as a function of two variables $m$ and $\alpha$. Then $r(m,\alpha)\ge r(m+1,\alpha)$ and $r(m,\alpha)\le r(m,\alpha+\epsilon)$ for any $\epsilon>0$.
\end{proposition}

\noindent\textbf{Remark.} We conjecture a stronger version of Proposition \ref{cou} stating the strict inequalities 
\[r(m,\alpha)> r(m+1,\alpha),\quad r(m,\alpha)< r(m,\alpha+\epsilon).\]

\begin{figure}[h!]
    \centering
    \includegraphics[scale=0.8]{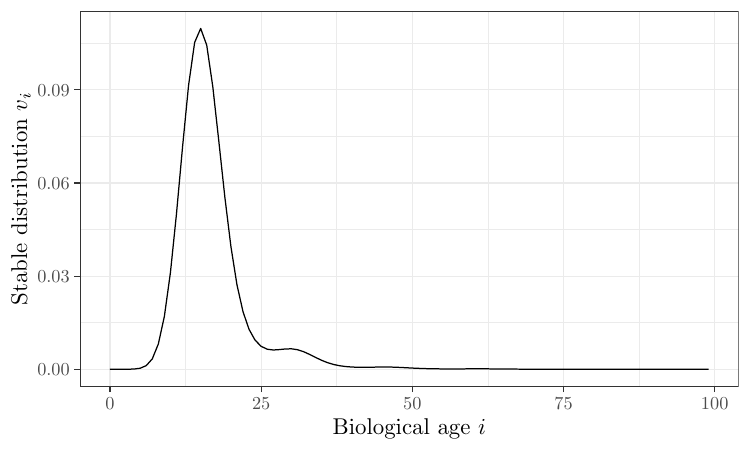}
    \caption{The stable distribution of biological ages $\boldsymbol{v}$ for $n = 100, m = 15,p = 0.5, \alpha = 2$.}
    \label{stab}
\end{figure}
\section{Average biological age}\label{PAA}
Given \eqref{tuv} with $r>1$, the left eigenvector $\boldsymbol{v}=(v_0,\ldots, v_{n-1})$ normalised such that 
\[v_0+\ldots+v_{n-1}=1,\]
represents the stable distribution over the type space for the individuals in the supercritical multitype Galton-Watson process, cf \cite{J}. Figure \ref{stab} illustrates the stable distribution corresponding to the case used to generate Figure \ref{ill}. This distribution should be inspected in comparison with the interval $19 \leq i \leq 63$ of rejuvenation states.

The related population mean value
\[a=\sum_{i=1}^{n-1} iv_i\]
can be treated as the average biological age across the cell population. The average biological age $a$ can be regarded as a measure of a population load related to biological aging: the lower the value of $a$, the fewer harmful proteins per cell are present in the population.

\begin{figure}[h!]
    \centering
    \includegraphics[width=0.9\textwidth]{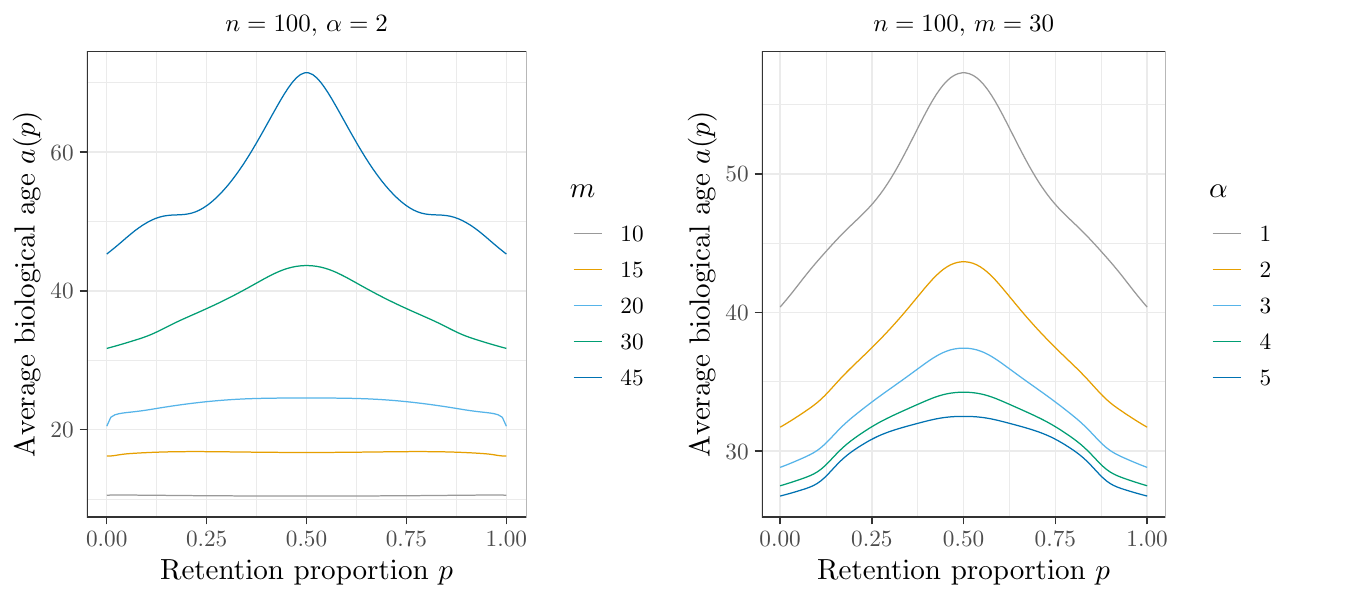}
    \caption{The average biological age $a = a(p)$ as a function of the retention proportion $p \in (0,1)$. Left panel: the five curves are generated for $n = 100, \alpha = 2$ and $m = 10, 15, 20, 30, 45$ in the ascending order of the curves. Right panel: the five curves are generated for $n = 100, m = 30$ and $\alpha = 1,2,3,4,5$ in the ascending order of the curves.}
    \label{mage}
\end{figure}

Figure \ref{mage} displays the average biological ages as a function of $p$ for various combinations of the model parameters $(n,m,\alpha)$. A noteworthy common feature across the majority of the depicted curves is the attainment of maximum population load at the equilibrium proportion $p$. However, it is crucial to emphasize that there are additional patterns, as demonstrated by the case $n=100$, $m=15$, $\alpha=2$.

In the left panel of Figure \ref{mage}, it can be observed that the average biological age $a$ increases in tandem with the constant inflow $m$. On the other hand, as depicted in the right panel, a higher division rate results in larger values of $a$. The last display supports the following statement.

\begin{proposition}\label{aou}
Given \eqref{spe}, \eqref{spm}, for a fixed triplet $(n,m,p)$, consider the average biological age $a=a(\alpha)$ as a function of the division rate $\alpha$. Then $a(\alpha)> a(\alpha+\epsilon)$ for any $\epsilon>0$.
\end{proposition}

\begin{figure}[h!]
\centering
\includegraphics[width=0.9\textwidth]{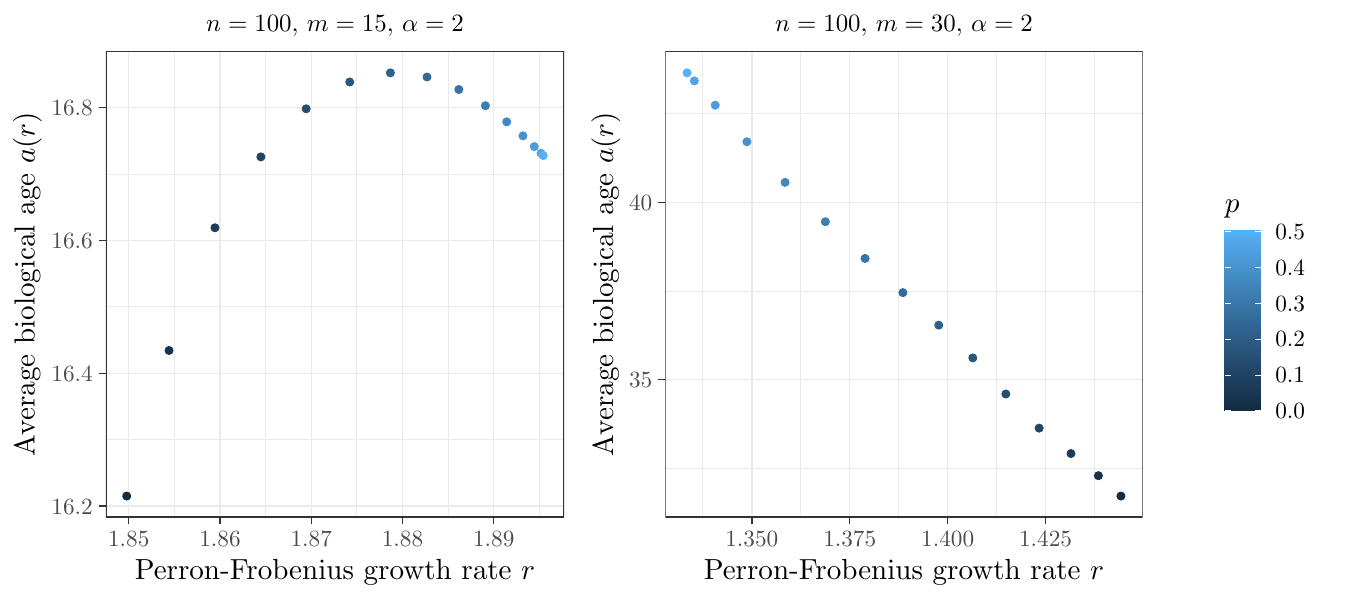}
\caption{The relation between the Perron-Frobenius growth rate $r=r(p)$ on the $x$-axis, and the average biological age $a=a(p)$ on the $y$-axis, across the different values of the retention proportion $0.025\le p\le 0.5$. Here the parameters $n=100$,  $\alpha=2$ are fixed. The two panels represents different values for $m$: on the left $m=30$ and on the right $m=50$.  }
 \label{fig2}
\end{figure}

\noindent\textbf{Remark}. Propositions \ref{cou} and \ref{aou} together with the right panel of Figure \ref{fig2} may give a false impression that a larger growth rate $r$ should always correspond to a smaller population load $a$. However, the left panel of Figure \ref{fig2} illustrates that the relationship between  $r$ and $a$ is more intricate.

Each panel of Figure \ref{fig2} hints at a line connecting the points on the plane, where each point $(x,y)=(r,a)$ represents the Perron-Frobenius growth rate and the average biological age corresponding to a certain set of parameters $(n,m,p,\alpha)$. Here the vector $(n,m,\alpha)$ is fixed, and the retention proportion $p$ takes different values across a broad range.
The more complex relationship depicted in the left panel reflects the specific pattern identified in the left panel of Figure \ref{mage} for the model parameters $n=100$, $m=15$, $\alpha=2$.

\section{Proofs}\label{proof}
\subsection*{Proof of Proposition \ref{MC}}

It is a straightforward exercise to obtain the stated expressions for the transition probabilities $p_{ij}$ for the Markov chain $\xi$ from the rule \eqref{trans}.
Observe that according to the stated formulas for $Q_{ij}$ and $R_{ij}=R_{ij}(p)$, 
\[\sum_{j=0}^nQ_{ij}=1-\sum_{j=i}^\infty q_{j-i}b_j,\quad 0\le j<n,\]
and, writing $g_{ij}$ instead of $g_{ij}(p)$, we obtain
\begin{align*}
    \sum_{j=0}^{n}R_{ij}&=\sum_{j=0}^{n-1}\sum_{k\ge i}q_{k-i}b_kg_{kj}+\sum_{k\ge n}q_{k-i}b_k\sum_{j=n}^{k}g_{kj}
=\sum_{k=i}^\infty q_{k-i}b_k.
\end{align*}
Thus, we have verified that 
$$\sum_{j=0}^{n}p_{ij}=\sum_{j=0}^nQ_{ij}+\sum_{j=0}^{n}R_{ij}=1,\quad 0\le i\le n.$$

\subsection*{Proof of Proposition \ref{coro}}
If $b_i=0$ for $i\ge n$, then by Proposition \ref{MC}, 
\begin{align*}
    Q_{ij}&=q_{j-i}(1-b_j),\quad R_{ij}=\sum_{k= j}^{n-1}q_{k-i}b_kg_{kj},\quad 0\le i\le j<n,\\
Q_{ij}&=0,\quad R_{ij}=\sum_{k= i}^{n-1}q_{k-i}b_{k}g_{kj},  \quad   0\le j<i<n,\\
Q_{in}&=\sum_{k\ge n}q_{k-i},\quad  R_{in}=0,\quad 0\le i<n.
\end{align*}
To compute $h_i$, observe first that 
\[\sum_{j=1}^{n}jQ_{ij}=\sum_{j=i}^{n-1}jq_{j-i}(1-b_j)+n\sum_{k\ge n}q_{k-i}=\rE((i+\tau)\wedge n -(i+\tau)b_{i+\tau}).\]
On the other hand, since $\sum_{j=1}^{k}j g_{kj}=kp$,
we have
\begin{align*}
    \sum_{j=1}^{n}jR_{ij}&=\sum_{j=1}^{n-1}\sum_{k=i}^{n-1}q_{k-i}b_kjg_{kj}=\sum_{k=i}^{n-1}q_{k-i}b_k\sum_{j=1}^{n-1}jg_{kj}\\
&=p\sum_{k=i}^{n-1}kq_{k-i}b_k=p\rE((i+\tau)b_{i+\tau}).
\end{align*}
Putting together the last two relations we arrive at \eqref{dav}, which in turn implies that the inequality $d_i<0$
is equivalent to  \eqref{day}. 

\subsection*{Proof of Proposition \ref{prog}}

Given \eqref{spe} and \eqref{spm}, relation \eqref{dav} entails the stated relations for $d_i$. For $0\le i<n-m$, it follows
\[n^{-1}d_i=n^{-1}m-(1-p)(y_i-y_i^{1+\alpha}),\]
where $y_i=\tfrac{i+m}{n}$.
Thus the inequality $d_i<0$ holds if and only if  $0\le i< n-m$ and
\[y_i-y_i^{1+\alpha}<\tfrac{m}{(1-p)n}.\]
This observation and the convexity of the function $y-y^{1+\alpha}$ over the interval $0<y<1$ yields the main statement of Proposition \ref{prog}.

\subsection*{Proof of Proposition \ref{cop}}

Let $0\le i<n-m$. It suffices to prove the following three inequalities
\begin{align}
 h_i(p,m,\alpha)&<  h_i(p+\epsilon,m,\alpha),\quad 0\le p<p+\epsilon\le1,\label{h1}\\
 h_i(p,m,\alpha)&< h_i(p,m+1,\alpha),\label{h2}\\
 h_i(p,m,\alpha)&> h_i(p,m,\alpha+\epsilon),\quad 1\le \alpha<\alpha+\epsilon. \label{h3}
\end{align}
Each of these inequalities is established using the coupling method \cite{Li}.

Turning to the proof of the first inequality, consider two versions of the Markov chain in question with different retention proportions: a version $\xi$ with $\xi_0=i$ and the retention proportions $p$, alongside a version $\xi'$ with $\xi'_0=i$ and the retention proportion $p'=p+\epsilon$. We arrange a coupling of these two versions in such a way that $\rP(\chi_1=\chi'_1)=1$, and observe that in the case of no division $\chi_1=\chi'_1=0$, we get $\rP(\xi_1<\xi'_1)=1$, while in the case of coupled divisions $\chi_1=\chi'_1=1$, the allocation of the harmful proteins according to the binomial scheme with $\xi_1\sim$ Bin$(i+m,p)$ and $\xi'_1\sim$ Bin$(i+m,p')$, may be arranged such that 
\begin{align}
\rP(\xi_1\le\xi'_1)&=1\text{ and }\rP(\xi_1<\xi'_1)>0,\label{one} 
\end{align} 
entailing \eqref{h1}. Inequality \eqref{h2} is established similarly.

To verify \eqref{h3}, we consider two versions of the Markov chain $\xi$ characterised by two different division rates: a version $\xi$ with $\xi_0=i$ and the division rate  $\alpha$, alongside a version $\xi'$ with $\xi'_0=i$ and the division rate $\alpha'=\alpha+\epsilon$. In this case, it is possible to arrange for
\begin{align}
\rP(\chi_1\ge\chi'_1)&=1\text{ and }\rP(\chi_1>\chi'_1)>0. \label{two} 
\end{align} 
If $\chi_1=\chi'_1=1$, then we can put $\xi_1=\xi'_1$. On the other hand, taking into account the outcome  $\chi'_1=0$, we can ensure that relation \eqref{one} holds. The relation \eqref{one} in turn implies \eqref{h3}.

\subsection*{Proof of Proposition \ref{coup}}

The four stated inequalities are proven using coupling. Here we prove the first inequality, while the other three are proven similarly to Proposition \ref{cop}.

Suppose $0\le i<j\le n-m$ and consider two coupled versions of the Markov chain $\xi$: $\xi$ starting from $\xi_0=i$ and $\xi'$ starting from $\xi'_0=j$. Due to the Markov property, it suffices to show that it is possible to organize a common probability space in such a way that \eqref{one} holds together with \eqref{two}.

It is obvious how to satisfy \eqref{two} using the strict monotonicity of the sequence \eqref{spe}. Furthermore, if $\chi'_1=0$, then obviously, $\rP(\xi_1<\xi'_1)=1$. Finally, if $\chi_1=\chi'_1=1$, then $\xi_1\sim$ Bin$(i+m,p)$ and $\xi'_1\sim$ Bin$(j+m,p)$, so that with $i< j$, it is easy to couple these two binomial distributions to ensure \eqref{one}.

\subsection*{Proof of Proposition \ref{linn}}
    In the $(m,\alpha)$-case, matrix \eqref{ma} has a decomposable form, 
\[
\mathbb{M} = \begin{pmatrix} \mathbb{\hat M}& \mathbb{M^*} \\ 0& 0 \end{pmatrix}
\]
(similar to the multitype Galton-Watson process of \cite{JS}, where the type of an individual is its chronological age). Since all elements of
$\mathbb{\hat M}$ 
are positive, the matrix $\mathbb{\hat M}$ has a Perron-Frobenius eigenvalue $r$, implying
\[\mathbb {\hat M}^t\sim r^t \boldsymbol{\hat u}^\intercal \boldsymbol{\hat v},\quad t\to\infty.\]
This observation yields the main statement of Proposition \ref{linn} since
\[
\mathbb{M}^t = \begin{pmatrix} \mathbb{\hat M}^t& \mathbb{\hat M}^{t-1}\mathbb{M^*} \\ 0& 0 \end{pmatrix}\sim r^t\begin{pmatrix} \boldsymbol{\hat u}^\intercal \boldsymbol{\hat v}& r^{-1}\boldsymbol{\hat u}^\intercal \boldsymbol{\hat v}\mathbb{M^*} \\ 0& 0 \end{pmatrix}=r^t \boldsymbol{u}^\intercal \boldsymbol{v}.
\]

\subsection*{Proof of Proposition \ref{cou}}
Even this proof is based on a coupling argument. To prove the inequality $r(m,\alpha)\ge r(m+1,\alpha)$, it suffices to prove that we can construct two $n$-type Galton-Watson processes $\boldsymbol{Z}_t$, driven by parameters $(n,m,p,\alpha)$, and $\boldsymbol{Z}'_t$, driven by parameters $(n,m+1,p,\alpha)$, each stemming from a single individual of the same type, and satisfying the inequality 
\begin{equation}\label{inz}
Z_t\ge Z'_t\text{ for all }t\ge0, 
\end{equation}
involving the total population sizes
\[Z_t=Z_t(0)+\ldots+Z_t(n-1),\quad Z'_t=Z'_t(0)+\ldots+Z'_t(n-1).\]
To this end, we compare two mother cell lineages: $\xi$, with parameters $(m,p,\alpha)$, and $\xi'$, with parameters $(m+1,p,\alpha)$. To establish \eqref{inz}, it is enough to verify that the two Markov chains can be constructed on the same probability space in such a way that conditionally on  $\xi_t=i\le j=\xi'_t$, 
\begin{equation}\label{p7}
\rP(\xi_{t+1}\le \xi'_{t+1})=1,\quad \rP( \eta_{t+1}\le \eta'_{t+1})=1.
\end{equation}
These relations ensure that all cell lineages in the process with parameters $(n,m,p,\alpha)$ are located further apart from the senescent state compared to the process with parameters $(n,m+1,p,\alpha)$.

The first equality in \eqref{p7} follows from the proof of Proposition \ref{coup}. To see that both equalities in \eqref{p7} hold, observe that with \eqref{spe}, the   $(n,m,p,\alpha)$-process dominates by the number of divisions over the $(n,m+1,p,\alpha)$-process. Furthermore, with more harmful proteins in the mother cell in the division event, one can arrange the allocation of proteins between the mother and the daughter cells in such a way that both mother and daughter cells receive a larger or equal number of proteins, entailing both relations in \eqref{p7}. 

The inequality $r(m,\alpha)\le r(m,\alpha+\epsilon)$ is proven similarly.

\subsection*{Proof of Proposition \ref{aou}}
This proposition can be readily demonstrated through a coupling argument, leveraging the intuitive feature of the model, which suggests that with a constant inflow of harmful proteins, each division reduces the number of harmful proteins per cell.

\section*{Acknowledgements} This research was partially supported by the Swedish Research Council (VR2023-04319) to MC. The funders had no role in study design, data collection and analysis, decision to publish, or preparation of the manuscript.

\end{document}